\documentclass{amsart}
\usepackage{amsopn}
\usepackage{mathtools}
\usepackage{mathabx}
\usepackage{amsmath}
\usepackage{amsthm}
\usepackage{amssymb, amscd}
\usepackage{physics}
\usepackage{tensor}
\usepackage{multirow}
\usepackage{graphicx, graphics, epsfig}
\usepackage{faktor} 
\usepackage{enumerate}
\usepackage{tikz}
\usepackage{hyperref}
\usepackage{color} 

\usepackage{amsfonts}
\usepackage{graphicx}
\usepackage{verbatim}
\immediate\write18{texcount -tex -sum  \jobname.tex > \jobname.wordcount.tex}

\DeclareMathOperator{\vol}{vol}

\DeclareMathOperator{\Span}{span}
\DeclareMathOperator{\Lie}{Lie}

\DeclareMathOperator{\DV}{div}
\DeclareMathOperator{\End}{End}

\newcommand{\nc}{\newcommand}
\nc{\SO}{\mathrm{SO}}
\nc{\Spin}{\mathrm{Spin}}
\nc{\Sp}{\mathrm{Sp}}
\nc{\Sl}{\mathrm{SL}}
\nc{\SU}{\mathrm{SU}}
\nc{\Or}{\mathrm{O}}
\nc{\U}{\mathrm{U}}
\nc{\GL}{\mathrm{GL}}
\nc{\Se}{\mathrm{S}}
\nc{\Cl}{\mathrm{Cl}} 
\nc{\Pin}{\mathrm{Pin}}
\nc{\pts}{\cdots\hspace{-0.1em}}
\nc{\nl}{\phantom{.} \newline}

\theoremstyle{plain}
\newtheorem{theorem}{Theorem}[section]
\newtheorem{proposition}[theorem]{Proposition}
\newtheorem{corollary}[theorem]{Corollary}
\newtheorem{lemma}[theorem]{Lemma}

\theoremstyle{definition}

\theoremstyle{remark}
\newtheorem{remark}[theorem]{Remark}

\title{New examples of $\mathrm{G}_2$-structures with divergence-free torsion}

\author{Agustín Garrone}  

\address{FAMAF, Universidad Nacional de Córdoba, Medina Allende s/n, ciudad universitaria (X5000HUA), Córdoba, Argentina.}
\email{agustin.garrone@mi.unc.edu.ar} 

\date{\today}

\begin{document}

\maketitle

\pagenumbering{arabic} 
\pagestyle{plain}

\begin{abstract}
    Interest in Riemannian manifolds with holonomy equal to the exceptional Lie group $\mathrm{G}_2$ have spurred extensive research in geometric flows of $\mathrm{G}_2$-structures defined on seven-dimensional manifolds in recent years. Among many possible geometric flows, the so-called \textit{isometric flow} has the distinctive feature of preserving the underlying metric induced by that $\mathrm{G}_2$-structure, so it can be used to evolve a $\mathrm{G}_2$-structure to one with the smallest possible torsion in a given metric class. This flow is built upon the divergence of the full torsion tensor of the flowing $\mathrm{G}_2$-structures in such a way that its critical points are precisely $\mathrm{G}_2$-structures with divergence-free torsion. In this article we study three large families of pairwise non-equivalent non-closed left-invariant $\mathrm{G}_2$-structures defined on simply connected solvable Lie groups previously studied in \cite{KL} and compute the divergence of their full torsion tensor, obtaining that it is identically zero in all cases.  
\end{abstract}


 
\section{Introduction}\label{intro}

\indent A $\mathrm{G}_2$-structure on a $7$-dimensional smooth manifold is a \textit{positive} smooth $3$-form $\varphi \in \Omega^3(M)$, by which we mean that for all $p \in M$ it can be written as
\begin{equation} \label{positiva intro}
    \varphi_p = e^{127} + e^{347} + e^{567} + e^{135} - e^{146} - e^{236} - e^{245}
\end{equation}
\noindent with respect to an ordered basis $\{e_1, \pts, e_7\}$ of $T_p M$. The existence of a $\mathrm{G}_2$-structure imposes several topological restrictions on $M$ beyond dimension, for example that $M$ be orientable and spin (\cite{La}, \cite{K1}, \cite{B1}, or \cite[Proposition 10.1.6]{J}). Notably, they are non-trivially related to the exceptional Lie group $\mathrm{G}_2$, the normed algebra of the octonions, generalizations of vector cross products, and even theoretical physics (most notably, $M$-theory); for a non-exhaustive review of these features, see \cite{Gr}, \cite{FG}, \cite[Section 2]{La}, \cite{B}, and \cite[Section 1]{G4}. 

\indent A striking result of $\mathrm{G}_2$-geometry is that a $\mathrm{G}_2$-structure naturally determines a Riemannian metric $g$ and a volume form $\vol$ in $M$; see Section \ref{preliminares}. In particular, they induce a Hodge star operator $\star:\Omega^k(M) \to \Omega^{7-k}(M)$, $0 \leq k \leq 7$. This induced Riemannian structure allows us to consider the Levi-Civita derivative $\nabla \varphi$ of $\varphi$, which is an interesting object in its own right. Of special interest is the \textit{torsion-free} case, which is characterized by the condition $\nabla \varphi = 0$. Torsion-free $\mathrm{G}_2$-structures are important for a variety of reasons, which include that the holonomy of $M$ with respect to $g$ is contained in $\mathrm{G}_2$ \textit{precisely} when $\nabla \varphi = 0$ (see \cite[Proposition 10.1.3]{J}), which in turn is relevant because of the famous Berger's classification theorem (see \cite{Be}). Historically, torsion-free $\mathrm{G}_2$-structures with holonomy equal to $\mathrm{G}_2$ were thought to be non-existent for decades until the first non-compact and compact examples were found respectively by R. Bryant in 1987 (see \cite{B1}) and D. Joyce in 1996 (see \cite{J}). Since then, geometric flows techniques have become a feature in the pursuit of new examples. Different geometric flows of $\mathrm{G}_2$-structures have been studied, notably the \textit{Laplacian flow of closed $\mathrm{G}_2$-structures} (i.e., such that $d \varphi = 0$) (see \cite{B2}) and the \textit{Laplacian coflow of coclosed $\mathrm{G}_2$-structures} (i.e., such that $d \star \varphi = 0$) (see \cite{KMMP}). The most relevant flow to us is the \textit{isometric flow}, first defined by S. Grigorian (see \cite{G2} and \cite{G3}) and independently by S. Dwivedi, P. Giannotis, and S. Karigiannis (see \cite{DGK}), which is governed by the equation
\begin{equation} \label{isometric flow intro}
    \begin{split}
    \frac{\partial \varphi(t)}{\partial t} &= \iota_{\DV T(t)}(\star \varphi(t)).
    \end{split}
\end{equation}
\noindent Several notions are introduced in equation \eqref{isometric flow intro}:
\begin{itemize}
    \item The map $\iota$, called \textit{contraction}, simply given by $\iota_X(\eta)(Y,Z, W) = \eta(X, Y, Z, W)$ for a $4$-form $\eta \in \Omega^4(M)$ and vector fields $X$, $Y$, $Z$, $W \in \mathfrak{X}(M)$.
    \item The \textit{full torsion tensor} $T$ of $\varphi$, which is the only smooth tensor field of type $(1,1)$ on $M$ such that (see Section \ref{preliminares}).
\begin{equation*}
    \nabla_X \varphi = \iota_{T(X)} (\star \varphi), \; \; X \in \mathfrak{X}(M).
\end{equation*}
    \item The \textit{divergence} of $T$, which is the smooth vector field $\DV T$ given by 
\begin{equation*}
    g(\DV T, E_j) = \sum_{i=1}^{7} (\nabla_{E_i}T) (E_i, E_j),
\end{equation*}
    \noindent where $\{E_i \, \vert \, 1 \leq i \leq 7\}$ is an arbitrary local orthonormal frame of $(M,g)$. 
\end{itemize}

\noindent It is clear that that $\mathrm{G}_2$-structures such that $\DV T = 0$ are critical points of equation \eqref{isometric flow intro}; such structures are said to have \textit{divergence-free full torsion tensor}, or \textit{divergence-free torsion} for short.  

\indent The main motivation behind this article is to find new examples of $\mathrm{G}_2$-structures with divergence-free full torsion tensor. Our research led us to left-invariant $\mathrm{G}_2$-structures defined on simply connected completely solvable Lie groups $G_{A,B,C}$, introduced in Section \ref{los ejemplos}, whose definition depends on the choice of a triple $A$, $B$, $C \in \mathfrak{sl}(4, \mathbb{R})$ of pairwise-commuting matrices. Recall from \cite{KL} that these groups admit lattices when $A$, $B$, $C$ are chosen to be linearly independent, so there are compact versions of the examples. Three broad families of triples are considered: the \textit{skew-symmetric case}, in which $A$, $B$, $C \in \mathfrak{sl}(4, \mathbb{R})$ are skew-symmetric; the \textit{diagonal case}, in which $A$, $B$, $C \in \mathfrak{sl}(4, \mathbb{R})$ are diagonal (which is proven to be equivalent to the \textit{symmetric case}, see Proposition \ref{prop: A,B,C sim compt}); and the \textit{antidiagonal case}, in which $A$, $B$, $C \in \mathfrak{sl}(4, \mathbb{R})$ are antidiagonal (see equation \eqref{eq: definicion antidiagonal}). These three cases are non-equivalent under mild restrictions on $A$, $B$, $C$, as computations regarding the Ricci operator of the underlying Riemannian structure of $(\mathfrak{g}_{A,B,C}, \varphi)$ readily show (see Proposition \ref{Ric in g_ABC}). The main result can be summarized as follows:
\begin{theorem}
     Left-invariant $\mathrm{G}_2$-structures given by equation \eqref{positiva intro} defined on certain simply connected completely solvable Lie groups $G_{A,B,C}$ (see Section \ref{los ejemplos}) with $A$, $B$, $C \in \mathfrak{sl}(4, \mathbb{R})$ pairwise-commuting matrices that are all skew-symmetric, all symmetric, or all antidiagonal have a divergence-free full torsion tensor $T$. Furthermore, those $\mathrm{G}_2$-structures are non-closed and pairwise non-equivalent under mild restrictions on $A$, $B$, and $C$.
\end{theorem}

\indent It seems natural to study the behavior of these $\mathrm{G}_2$-structures under an appropriate flow, mirroring the work in \cite[Sections 4.6 and 4.7]{KL} for the Laplacian coflow of a class of coclosed $\mathrm{G}_2$-structures defined over $\mathfrak{g}_{A,B,C}$ and a modification of such flow. Another natural question is wheter these $\mathrm{G}_2$-structures are stable or unstable critical points of the \textit{energy functional} considered in \cite{DGK} and \cite{G2}. We leave these questions open for future works. 


\section{\texorpdfstring{$\mathrm{G}_2$}{G2}-structures and their geometry} \label{preliminares}

\subsection{Definitions}

\indent A smooth $3$-form $\varphi$ defined on a $7$-dimensional differentiable manifold $M$ is called a $\mathrm{G}_2$-structure if it is \textit{positive}; that is, if at every point $p \in M$ it can be written as
\begin{equation} \label{positiva}
    \varphi_p = e^{127} + e^{347} + e^{567} + e^{135} - e^{146} - e^{236} - e^{245}
\end{equation}   
\noindent with respect to an ordered basis $\{e_1, \pts, e_7\}$ of $T_p M$. Two $\mathrm{G}_2$-structures $(M,\varphi)$ and $(M', \varphi')$ are said to be \textit{equivalent} if there is a diffeomorphism $f:M \to M'$ such that $\varphi = f^* \varphi'$; such a map is called an \textit{equivalence}. 

\indent One of the most astounding and well-known facts about $\mathrm{G}_2$-structures is that they canonically induce a Riemannian metric $g$ and a volume form $\vol$ by
\begin{equation} \label{g y vol en G2}
    g(X,Y) \vol = \tfrac{1}{6} \iota_X(\varphi) \wedge \iota_Y(\varphi) \wedge \varphi, \quad X, Y \in \mathfrak{X}(M);
\end{equation} 
\noindent here, the map $\iota$ is simply defined as $\iota_X(\varphi)(W,Z) := \varphi(X,W,Z)$ for $W$, $Z \in \mathfrak{X}(M)$. This follows from straighforward computation in coordinates, as is oulined for instance in \cite[Section 2, Theorem 1]{B1}. $g$ and $\vol$ induce a Hodge star operator  $\star:\Omega^k(M) \to \Omega^{7 - k}(M)$, $0 \leq k \leq 7$, on $M$; in particular, it allows us to define the dual $4$-form of $\varphi$, $\psi := \star \varphi$, which can readily be shown to satisfy
\begin{equation*}
    \psi_p = e^{3456} + e^{1256} + e^{1234} - e^{2467} + e^{2357} + e^{1457} + e^{1367}
\end{equation*}
\noindent with respect to the ordered basis $\{e_1, \pts, e_7\}$. Note that our definition \eqref{positiva} of $\varphi$ differs from other more commonly used alternatives; we mention that the one used in \cite{B1} and \cite{B2} is fact equivalent to ours, and that the one used in \cite{K2} and \cite{K1} induces the opposite orientation to ours.

\indent Some of the most important aspects of the geometry that can be done in presence of $\mathrm{G}_2$-structures are encoded in the \textit{torsion forms}, which are the unique differential forms $\tau_i \in \Omega^i(M)$ with $i = 0$, $1$, $2$, $3$ such that
\begin{equation} \label{eq: phi phis}
    d \varphi = \tau_0 \psi + 3 \tau_1 \wedge \varphi + \star \tau_3, \quad d \psi = 4 \tau_1 \wedge \psi - \star \tau_2;
\end{equation}
\noindent this is proven for example in \cite[Section 3.3, Proposition 1]{B2}. Equivalently, we have 
\begin{equation} \label{eq: tau taus}
    \begin{array}{c}
       \tau_0  = \tfrac{1}{7} \star (d \varphi \wedge \varphi), \quad  \tau_1  = - \tfrac{1}{12} \star (\star d \varphi \wedge \varphi), \\
        \tau_2  = - \star d \psi + 4 \star (\tau_1 \wedge \psi), \quad \tau_3 = \star d \varphi - \tau_0 \varphi - 3 \star (\tau_1 \wedge \varphi);
    \end{array}
\end{equation}
\noindent see \cite[Section 1.2]{MOV} or \cite[Section 2.1]{LMSES} for a more careful derivation of this equivalence. The vanishing of some torsion forms determines different kinds of $\mathrm{G}_2$-structures, for example:
\begin{itemize}
    \item $d \varphi = 0$ if and only if $\tau_0 = \tau_1 = \tau_3 = 0$, in which case we say $\varphi$ is \textit{closed}. 
    \item $d \psi = 0$ if and only if $\tau_1 = \tau_2 = 0$, in which case we say $\varphi$ is \textit{coclosed}. 
    \item $d \varphi = 0$ and $d\psi = 0$ if and only if $\tau_i = 0$ for all $i = 0$, $1$, $2$, $3$, in which case we say $\varphi$ is \textit{torsion-free}. 
\end{itemize} 

\indent Recall that the \textit{torsion} of a $\mathrm{G}_2$-structure is $\nabla \varphi$, the covariant derivative of $\varphi$ with respect to the metric $g$ induced by $\varphi$. It can be shown that there is a tensor field $T \in \mathfrak{T}^{(1,1)}(M)$ such that
\begin{equation*} 
    \nabla_X \varphi = \iota_{T(X)} (\psi), \quad X \in \mathfrak{X}(M);
\end{equation*}
\noindent we refer to \cite[Theorem 4.44]{K1} or \cite[Theorem 2.27]{K2} for a proof. We call $T$ the \textit{full torsion tensor} of $\varphi$. As the name suggests, $T$ is directly related to the torsion forms $\tau_i \in \Omega^i(M)$, $i = 0$, $1$, $2$, $3$, as it can be shown that
\begin{equation} \label{eq: full torsion tensor form}
     T(X,Y) =  \tfrac{1}{4} \tau_0 g(X,Y) - \iota_{\tau_1}(\varphi)(X,Y) -  \tfrac{1}{2} \tau_2(X,Y) - \tau_{27}(X,Y), \quad X, Y \in \mathfrak{X}(M),
\end{equation}
\noindent where $\iota_{\tau_1}(\varphi)$ is the contraction of $\varphi$ by the smooth vector field $\widetilde{\tau}_1$ uniquely determined by the condition 
\begin{equation*}
    g(\widetilde{\tau}_1, X) = \tau_1(X), \quad X, \in \mathfrak{X}(M),
\end{equation*}
\noindent which is called $\tau_1$ by abuse of language, and $\tau_{27}$ is the symmetric 2-tensor given by
\begin{equation} \label{tau_27}
    \tau_{27}(X,Y) := \star(\iota_X(\varphi) \wedge \iota_Y(\varphi) \wedge \tau_3), \quad X, Y \in \mathfrak{X}(M);
\end{equation}
\noindent a proof can be found in \cite[Theorem 2.27]{K2}. We see that a $\mathrm{G}_2$-structure is torsion-free if and only if $T = 0$, which is also equivalent to $\nabla \varphi = 0$. This is often called Fernández-Gray theorem; see \cite{FG}.  

\indent Interest in torsion-free $\mathrm{G}_2$-structures is widespread, as these structures have many desirable properties, among them having holonomy contained in the exceptional Lie group $\mathrm{G}_2$ and being Ricci-flat (see \cite[Propositions 10.1.3 and 10.1.5]{J}), and are of historical relevance, as the famous Berger's classification theorem illustrates (see \cite{Be}). Geometric flows techniques have been introduced in this area of research, most notably in \cite{B2} and in \cite{KMMP}. Of special relevance to us is the \textit{isometric flow}, first defined by S. Grigorian (see \cite{G2} and \cite{G3}) and independently by S. Dwivedi, P. Giannotis, and S. Karigiannis (see \cite{DGK}), which is governed by the equation
\begin{equation} \label{isometric flow}
    \begin{split}
    \frac{\partial \varphi(t)}{\partial t} &= \iota_{\DV T(t)}(\psi(t));
    \end{split}
\end{equation}
\noindent here the \textit{divergence} $\DV T$ of the full torsion tensor $T$ is the vector field defined as 
\begin{equation} \label{eq: divergencia es nabla}
    g(\DV T, E_j) = \sum_{i=1}^{7} (\nabla_{E_i}T) (E_i, E_j), 
\end{equation}
\noindent where $\{E_i \, \vert \, 1 \leq i \leq 7\}$ is an arbitrary local orthonormal frame with respect to the induced metric $g$. It is clear that $\mathrm{G}_2$-structures such that $\DV T = 0$ are critical points of equation \eqref{isometric flow}; such structures are said to have \textit{divergence-free full torsion tensor}, or \textit{divergence-free torsion} for short. It is known from \cite[Theorem 4.3]{G2} that $\DV T = 0$ if $d \varphi = 0$ (see also Remark \ref{obs: closed are divergence-free}). The $\mathrm{G}_2$-structures considered in this article are non-closed in general. 

\subsection{Left-invariant structures}

\indent We will restrict our attention to left-invariant $\mathrm{G}_2$-structures defined on Lie groups $G$, for it is usual to come across simpler formulas in this setting. An important observation in this regard is that $\nabla_X Y \in \Lie(G)$ for all $X$, $Y \in \Lie(G)$, which in turn allows us to write the Koszul formula for the Levi-Civita connection as
\begin{equation} \label{eq: formula de Koszul linealizada}
        g(\nabla_X Y, Z) =  \tfrac{1}{2} \left( g([X,Y],Z) - g([Y,Z],X) + g([Z,X],Y) \right), \quad X, Y, Z \in \Lie(G).
\end{equation}
\noindent It is natural then to consider the map $U:\Lie(G) \times \Lie(G) \to \Lie(G)$ determined by
\begin{equation} \label{eq: U(X,Y)}
    g(U(X,Y),Z) =  \tfrac{1}{2} \left( g([Z,X],Y) - g([Y,Z],X) \right), \quad  X, Y, Z \in \Lie(G).
\end{equation}
\noindent Note in particular that $U(X,Y) = U(Y,X)$ for all $X$, $Y \in \Lie(G)$. In terms of this map, equation \eqref{eq: formula de Koszul linealizada} can be written more simply as
\begin{equation} \label{eq: nabla linealizado}
    \nabla_X Y =  \tfrac{[X,Y]}{2} + U(X,Y), \quad X, Y \in \Lie(G);
\end{equation}
\noindent equation \eqref{eq: nabla linealizado} will come in handy. It is worth mentioning that $U$ is identically zero when the metric $g$ is bi-invariant; that is, $g$ is both left and right-invariant.

\indent Turning again to $\mathrm{G}_2$-structures, we first recall that two left-invariant $\mathrm{G}_2$-structures $(G, \varphi)$ and $(G', \varphi')$ defined on Lie groups $G$ and $G'$ are called \textit{equivariantly equivalent} if there is a Lie group isomorphism $f:G \to G'$ such that $f^* \varphi' = \varphi$; such a map is called an \textit{equivariant equivalence}. We also recall that any $3$-form defined on the Lie algebra $\Lie(G)$ of a Lie group $G$ by
\begin{equation*}
    \varphi = e^{127} + e^{347} + e^{567} + e^{135} - e^{146} - e^{236} - e^{245} \in \Lambda^3(\Lie(G)^*),
\end{equation*}
\noindent with respect to some ordered basis $\{e_1, \pts, e_7\}$ on $\Lie(G)$ gives rise to a left-invariant $\mathrm{G}_2$-structure on $G$, typically denoted by the same name; moreover, $\varphi$ induces on $\Lie(G)$ an inner product $\langle \cdot, \cdot \rangle$ and a volume form $\vol$ on $\Lie(G)$ given by a formula akin to \eqref{g y vol en G2}, which can in turn be extended to corresponding left-invariant structures on the Lie group $G$. Similar considerations apply to other $\mathrm{G}_2$-related paraphernalia, most importantly the full torsion tensor $T$, whose divergence satisfies
\begin{equation} \label{eq: div alg lie}
    \langle \DV T,e_j \rangle = - \sum_{1 \leq i \leq 7} T(\nabla_{e_i} e_i, e_j) - \sum_{1 \leq i \leq 7} T(e_i, \nabla_{e_i} e_j)
\end{equation}
\noindent when restricted to $\Lie(G)$, where $\{e_1,\pts, e_7\}$ is an ordered orthonormal basis of $\Lie(G)$. This prompts us to work entirely at the Lie algebra level. 

\section{\texorpdfstring{$\mathrm{G}_2$}{G2}-structures over the Lie algebras \texorpdfstring{$\mathfrak{g}_{A,B,C}$}{g\_\{A, B, C\}}} \label{los ejemplos}

\subsection{Definitions and remarks}

\indent Let $A$, $B$, $C \in \mathfrak{sl}(4, \mathbb{R})$ be pairwise-commuting. Over a real seven-dimensional vector space with basis $\{e_1, \pts, e_7\}$, which we denote $\mathfrak{g}_{A,B,C}$, we put
\begin{equation*}
    \mathfrak{a} := \Span\{e_1, e_2, e_7\}, \quad \mathfrak{n} := \Span\{e_3, e_4, e_5, e_6\},
\end{equation*}
\noindent and define a bilinear, anticommutative product $[\cdot, \cdot]: \mathfrak{g}_{A,B,C} \times \mathfrak{g}_{A,B,C} \to \mathfrak{g}_{A,B,C}$ by setting $[\mathfrak{a}, \mathfrak{a}] = 0$, $[\mathfrak{n}, \mathfrak{n}] = 0$, and
\begin{equation} \label{eq: commutators g_A,B,C}
    [e_7, v] = Av, \quad  [e_1, v] = Bv, \quad [e_2, v] = Cv
\end{equation}    
\noindent for all $v \in \mathfrak{n}$, so as to make $(\mathfrak{g}_{A,B,C}, [\cdot, \cdot])$ a real $7$-dimensional Lie algebra. 
    
\indent These algebras were first defined in \cite[Section 4]{KL}, and they constitute a sort-of generalization of almost abelian Lie algebras. Note that the Jacobi condition is actually equivalent to the fact that $A$, $B$, $C$ are pairwise-commuting. By definition, $\mathfrak{a} = \Span \{e_1, e_2, e_7 \}$ is an abelian subalgebra and $\mathfrak{n} = \Span \{e_3, e_4, e_5, e_6 \}$ is an abelian ideal of $\mathfrak{g}_{A,B,C}$; furthermore, it is readily verified that $\mathfrak{n}$ is the nilradical of $\mathfrak{g}_{A,B,C}$ only if no linear combination of $A$, $B$, $C$ is a nilpotent matrix. For all choices of $A$, $B$, $C$ we have that $\mathfrak{g}_{A,B,C}$ is a completely solvable Lie algebra and, as a consequence of tracelessness, that $\mathfrak{g}_{A,B,C}$ is unimodular. Under the further restriction that $A$, $B$, $C \in \mathfrak{sl}(4, \mathbb{R})$ are linearly independent and that they simultaneously diagonalize over $\mathbb{R}$, in which case we call the triple \textit{compatible}, all Lie algebras $\mathfrak{g}_{A,B,C}$ are isomorphic to each other, and the corresponding simply connected Lie groups $G_{A,B,C}$ are isomorphic to $G_J$, the only completely solvable unimodular Lie group in the classification appearing in \cite{LN}. We refer to \cite[Section 4]{KL} for more details on this topic. We will not be dealing with compatible triples in this article except for a brief remark in Section \ref{skew case}.

\indent We consider over $\mathfrak{g}_{A,B,C}$ the $\mathrm{G}_2$-structure defined with respect to the basis $\{e_1, \pts, e_7\}$ as 
\begin{equation} \label{eq: G2-estructura}
        \varphi := e^{127} + e^{347} + e^{567} + e^{135} - e^{146} - e^{236} - e^{245} \in \Lambda^3(\mathfrak{g}_{A,B,C}^*),
\end{equation}
\noindent which, as we know, induces an inner product $\langle \cdot, \cdot \rangle$ and a volume form $\vol$ on $\mathfrak{g}_{A,B,C}$ according to
\begin{equation*}
    g(X,Y) \vol =  \tfrac{1}{6} \iota_X(\varphi) \wedge \iota_Y(\varphi) \wedge \varphi, \quad X, Y \in \mathfrak{g}_{A,B,C}.
\end{equation*}
\noindent These definitions ensure that $\{e_1, \pts, e_7\}$ is orthonormal with respect to $\langle \cdot, \cdot \rangle$, as well as $\mathfrak{a}$ and $\mathfrak{n}$ are orthogonal complements to each other. Accordingly, this structure induces a Hodge star operator \\ $\star:\Lambda^k(\mathfrak{g}_{A,B,C}^*) \to \Lambda^{7-k}(\mathfrak{g}_{A,B,C}^*)$, which in turn allows us to define the dual $4$-form $\psi := \star \varphi$ as
\begin{equation*}
    \psi := e^{3456} + e^{1256} + e^{1234} - e^{2467} + e^{2357} + e^{1457} + e^{1367}. 
\end{equation*}
\indent In \cite[Section 4.1]{KL} it is established that the notion of equivariant equivalence among left-invariant $\mathrm{G}_2$-structures on simply connected Lie groups $G_{A,B,C}$ with Lie algebras of the form $\mathfrak{g}_{A,B,C}$ correspond with Lie algebras isomorphisms $h: \mathfrak{g}_{A,B,C} \to \mathfrak{g}_{A',B',C'}$ such that $h \in \mathrm{G}_2$, where
\begin{align*}
    \mathrm{G}_2 = \{h \in \GL(7, \mathbb{R}) \, \vert \, h^* \varphi = \varphi\}. 
\end{align*}
\noindent This result allows us not only to work entirely at the Lie algebra level but also with a fixed $\varphi$. With the help of this result we will be able to show that some choices of pairwise-commuting matrices $A$, $B$, $C \in \mathfrak{sl}(4, \mathbb{R})$ lead to non-equivalent examples of left-invariant $\mathrm{G}_2$-structures with divergence-free full torsion tensor defined on the simply connected Lie group $G_{A,B,C}$ and whose $\mathrm{G}_2$-structure is induced by equation \eqref{eq: G2-estructura}. 

\subsection{Formulas for the torsion forms} \label{formulas for the torsion}

\indent In order to simplify future formulas and calculations, we introduce the $2$-forms on $\mathfrak{g}_{A,B,C}$ given by
\begin{align} 
    \omega_7 := e^{34} + e^{56} , \quad
    \omega_1 := e^{35} - e^{46} , \quad
    \omega_2 := - e^{36} - e^{45}, \label{eq: omegas formulas} \\ 
    \overline{\omega}_7 := e^{34} - e^{56}, \quad
    \overline{\omega}_1 := e^{35} + e^{46} , \quad
    \overline{\omega}_2 := -e^{36} + e^{45}. \label{eq: omegas barra formulas}
\end{align}
\noindent We interpret the $2$-forms given by equations \eqref{eq: omegas formulas} and \eqref{eq: omegas barra formulas} as elements of $\Lambda^2(\mathfrak{n}^*)$. We now list some immediate properties of these $2$-forms in the following result for future reference. We call $g_{\mathfrak{n}}$, $\vol_{\mathfrak{n}}$, and $\star_{\mathfrak{n}}$ the metric, volume form, and Hodge star operator induced on $\mathfrak{n}$ by the corresponding structures in $\mathfrak{g}_{A,B,C}$. 

\begin{lemma} \label{lemma: omegas}
    Let $\mathfrak{n} = \Span\{e_3, e_4, e_5, e_6\}$ be the nilpotent ideal of $\mathfrak{g}_{A,B,C}$, and let $\omega_7$, $\omega_1$, $\omega_2$, $\overline{\omega}_7$, $\overline{\omega}_1$, and $ \overline{\omega}_2 \in \Lambda^2(\mathfrak{n}^*)$ be given as in equations \eqref{eq: omegas formulas} and \eqref{eq: omegas barra formulas}. Then
    \begin{enumerate} [\rm (i)]
        \item \label{omegas 1} $\varphi = e^{127} + \omega_7 \wedge e^7 + \omega_1 \wedge e^1 + \omega_2 \wedge e^2$.
        \item $\psi = e^{3456} + \omega_7 \wedge e^{12} + \omega_1 \wedge e^{27} - \omega_2 \wedge e^{17}$.
        \item $\star_{\mathfrak{n}} \omega_i = \omega_i$ and $\star_{\mathfrak{n}} \overline{\omega}_j = - \overline{\omega}_j$ for $i,j = 1, 2, 7$.
        \item $\omega_i \wedge \omega_j = \omega_i \wedge \overline{\omega}_j = \overline{\omega}_i \wedge \overline{\omega}_j = 0$ for all $i \neq j$.
        \item the following equations hold:
\begin{align*}
        \omega_7 \wedge \omega_7 &= \omega_1 \wedge \omega_1 = \omega_2 \wedge \omega_2 = \phantom{+} 2 e^{3456}, \\ \overline{\omega}_7 \wedge \overline{\omega}_7 &= \overline{\omega}_1 \wedge \overline{\omega}_1 = \overline{\omega}_2 \wedge \overline{\omega}_2 = - 2 e^{3456}.
\end{align*}
        \item the following equations hold:
\begin{align*}
    e^{34} &=  \tfrac{\overline{\omega}_7 + \omega_7}{2}, \quad
    e^{35} =  \tfrac{\overline{\omega}_1 + \omega_1}{2}, \quad
    e^{36} = -  \tfrac{\overline{\omega}_2 + \omega_2}{2}, \quad \\
    e^{45} &=  \tfrac{\overline{\omega}_2 - \omega_2}{2}, \quad
    e^{46} =  \tfrac{\overline{\omega}_1 - \omega_1}{2}, \quad 
    e^{56} = -  \tfrac{\overline{\omega}_7 - \omega_7}{2}.
\end{align*}
        \item $\mathfrak{B} := \{\overline{\omega}_7, \overline{\omega}_1, \overline{\omega}_2, \omega_7, \omega_1, \omega_2\}$ is an orthogonal basis of $\Lambda^2(\mathfrak{n}^*)$ with respect to the metric induced by $g_{\mathfrak{n}}$, with every element having norm equal to $\sqrt{2}$.
    \end{enumerate}
\end{lemma} 

\indent Let $\theta:\mathfrak{sl}(4, \mathbb{R}) \to \End(\Lambda^2(\mathfrak{n}^*))$ denote the natural representation of $\mathfrak{sl}(4, \mathbb{R})$ on $\Lambda^2(\mathfrak{n}^*)$, given by
\begin{equation*} 
    (\theta(D)\eta)(\cdot, \cdot) := - \eta(D \cdot, \cdot) - \eta(\cdot, D \cdot).
\end{equation*} 
\noindent It is clear from the definition that $\theta(M) e^i = - \sum_{j=3}^6 m_{ij} e^j$, where $M \in \mathfrak{sl}(4, \mathbb{R})$ is a matrix with coefficients $[m_{ij}]$, $i,j = 3, 4, 5, 6$. Moreover,
\begin{align} 
    \theta(M) \omega_7 :=& - (m_{33} + m_{44}) e^{34} + (m_{63} - m_{45}) e^{35} - (m_{46} + m_{53}) e^{36} \label{eq: theta en las omegas 1} \\
        &+ (m_{64} + m_{35})e^{45} + (m_{36} - m_{54})e^{46} - (m_{55} + m_{66})e^{56}, \notag \\
    \theta(M) \omega_1 :=& - (m_{54} + m_{63})e^{34} - (m_{33} - m_{55})e^{35} + (m_{43} - m_{56})e^{36} \label{eq: theta en las omegas 2} \\
        &+ (m_{65} - m_{34})e^{45} + (m_{44} + m_{66})e^{46} + (m_{45} + m_{36})e^{56},  \notag \\
    \theta(M) \omega_2 :=& \quad \; (m_{64} - m_{53})e^{34} + (m_{43} + m_{65})e^{35} + (m_{33} + m_{66})e^{36} \label{eq: theta en las omegas 3} \\
        &+ (m_{44} + m_{55})e^{45} + (m_{34} + m_{56})e^{46} + (m_{35} - m_{46})e^{56}; \notag 
\end{align}
\noindent see \cite[Remark 2.8]{N}\footnote{The reader should be alerted that equations \eqref{eq: theta en las omegas 3} and \eqref{tau0 N2} differ from their counterparts in \cite{N}. We believe this is due to typos present in the reference cited, and that our versions are correct.}. We recall the following result.

\begin{proposition} \cite[Theorems 2.3 and 2.4]{N} \label{prop: formulas} 
    For each set of pairwise-commuting matrices $A$, $B$, $C \in \mathfrak{sl}(4, \mathbb{R})$, the following formulas hold in $(\mathfrak{g}_{A,B,C}, \varphi)$:
\begin{align*} 
    d\varphi &= (\theta(B)\omega_7 - \theta(A)\omega_1) \wedge e^{17} + (\theta(C)\omega_7-\theta(A)\omega_2)\wedge e^{27} + (\theta(B)\omega_2-\theta(C)\omega_1)\wedge e^{12}, \\ 
    \star d\varphi &= (\theta(B^{\intercal})\omega_7-\theta(A^{\intercal})\omega_1)\wedge e^2-(\theta(C^{\intercal})\omega_7-\theta(A^{\intercal})\omega_2)\wedge e^1 -(\theta(B^{\intercal})\omega_2-\theta(C^{\intercal})\omega_1)\wedge e^7, \\ 
    d\psi &= \left( \theta(A)\omega_7+\theta(B)\omega_1+\theta(C)\omega_2 \right)\wedge e^{127},   \\ 
    \star d \psi &= -\left(\theta(A^{\intercal})\omega_7+\theta(B^{\intercal})\omega_1+\theta(C^{\intercal})\omega_2\right).  
\end{align*}    
\end{proposition}

\begin{remark}
    The formulas in Proposition \ref{prop: formulas} show that the $\mathrm{G}_2$-structures $(\mathfrak{g}_{A,B,C},\varphi)$ are neither closed nor coclosed in general. In fact, $(\mathfrak{g}_{A,B,C},\varphi)$ is closed if and only if 
\begin{equation*}
    \theta(B)\omega_7 = \theta(A) \omega_1, \quad \theta(C) \omega_7 = \theta(A) \omega_2, \quad \theta(B) \omega_2 = \theta(C) \omega_1,
\end{equation*}    
    \noindent and $(\mathfrak{g}_{A,B,C},\varphi)$ is coclosed if and only if $\theta(A)\omega_7+\theta(B)\omega_1+\theta(C)\omega_2 = 0$. 
\end{remark}

\indent Proposition \ref{prop: formulas} can be used to give explicit formulas for the torsion forms $\tau_0$, $\tau_1$, $\tau_2$, and $\tau_3$ of $(\mathfrak{g}_{A,B,C}, \varphi)$ in terms of the action of $\theta$ via equations \eqref{eq: tau taus}. This has been done in \cite[Proposition 2.7]{N} with the help of equations \eqref{eq: theta en las omegas 1} to \eqref{eq: theta en las omegas 3}; we recall the ugly-looking result\footnote{There is no explicit formula for the torsion form $\tau_3$ in \cite[Proposition 2.7]{N}, although clear guidelines for obtaining it are given; we have computed $\tau_3$ following those with the help of Maple 2023 from Maplesoft (see Appendix \ref{section: maple 2023} for a brief description of the code used).}.

\begin{proposition} \cite[Proposition 2.7]{N} \label{prop: formulas N2}
    For each set of pairwise-commuting matrices $A = [a_{ij}]$, $B = [b_{ij}]$, $C = [c_{ij}]\in \mathfrak{sl}(4, \mathbb{R})$ for $i$, $j \in \{3, 4, 5, 6\}$, the following formulas hold in $(\mathfrak{g}_{A,B,C}, \varphi)$:
\begin{align}
    \tau_0 =& \phantom{+} \; \, \tfrac{2}{7} ( a_{34} - a_{43} + a_{56} - a_{65} + b_{35} - b_{53} + b_{64} - b_{46} + c_{54} - c_{45} + c_{63} - c_{36} ), \label{tau0 N2} \\
    \tau_1 =& - \tfrac{1}{12} ( a_{36} - a_{63} + a_{45} - a_{54} +  c_{56} - c_{65} + c_{34} - c_{43} ) e^1 \label{tau1 N2}  \\  
    &-  \tfrac{1}{12} ( a_{64} - a_{46} + a_{35} - a_{53} +  b_{43} - b_{34} + b_{65} - b_{56} ) e^2 \notag \\  
    &-  \tfrac{1}{12} (b_{63} - b_{36} + b_{54} - b_{45} + c_{46} - c_{64} + c_{53} - c_{35}) e^7, \notag
\end{align}     
\begin{align} 
    \tau_2 =&\phantom{+} \; \, \tfrac{1}{3} (b_{45} - b_{54} + b_{36} - b_{63} + c_{35} - c_{53} + c_{64} - c_{46}) e^{12} \label{tau2 N2} \\ 
        &+  \tfrac{1}{3} (a_{64} - a_{46} + a_{35} - a_{53} + b_{65} - b_{56} + b_{43} - b_{34}) e^{17} \notag \\ 
        &+  \tfrac{1}{3} ( a_{54} - a_{45} + a_{63} - a_{36} + c_{65} - c_{56} + c_{43} - c_{34} ) e^{27} \notag \\  
        &+  \tfrac{1}{3} (- 3 a_{33} - 3 a_{44} + 2 c_{46} - 2 c_{35} - 2 b_{45} - 2 b_{36} - c_{53} + c_{64} - b_{63} - b_{54}) e^{34} \notag \\  
        &+  \tfrac{1}{3} (- 2 a_{54} + 2 a_{36} + 2c_{56} + 2c_{34} + a_{63} - a_{45} + c_{65} + c_{43} - 3 b_{55} - 3 b_{33}) e^{35} \notag \\  
        &+  \tfrac{1}{3} (- 2a_{64} - 2a_{35} - 2 b_{65} + 2b_{34} - a_{46} - a_{53} - b_{56} + b_{43} - 3c_{44} - 3c_{55}) e^{36} \notag \\  
        &+  \tfrac{1}{3} ( \phantom{+} a_{64} + a_{35} + b_{65} - b_{34} + 2a_{46} + 2a_{53} + 2b_{56} - 2b_{43} + 3c_{55} + 3c_{44}) e^{45} \notag \\  
        &+  \tfrac{1}{3} (-a_{54} + a_{36} + c_{56} + c_{34} + 2a_{63} - 2a_{45} + 2c_{65} + 2c_{43} - 3b_{33} - 3b_{55}) e^{46} \notag \\  
        &+  \tfrac{1}{3} (\phantom{+} 3 a_{33} + 3 a_{44} - c_{46} + c_{35} + b_{45} + b_{36} + 2c_{53} - 2c_{64} + 2b_{63} + 2b_{54}) e^{56}, \notag \\
    \tau_3 =& - \tfrac{2}{7} \left(a_{56} +c_{54} +a_{34} -c_{36} -a_{65} +c_{63} -a_{43} -c_{45} +b_{64} +b_{35} -b_{53} -b_{46} \right) e^{127} \label{tau3 N2} \\
        &- \tfrac{1}{4} \left(-b_{65} -b_{43} +b_{56} +b_{34} -a_{64} +a_{53} -3 a_{46} +3 a_{35} -c_{33} -c_{44} \right) e^{134} \notag \\
        &+ \tfrac{1}{7} \left(5 a_{56} +5 c_{54} +5 a_{34} -5 c_{36} +2 a_{65} -2 c_{63} +2 a_{43} +2 c_{45} -2 b_{64} -2 b_{35} +2 b_{53} +2 b_{46} \right) e^{135} \notag \\
        &- \tfrac{1}{4} \left(-b_{63} +b_{45} -b_{54} +b_{36} -c_{53} -c_{46} -3 c_{64} -3 c_{35} +a_{55} +a_{44} \right) e^{136} \notag \\
        &+ \tfrac{1}{4} \left(b_{63} -b_{45} +b_{54} -b_{36} -3 c_{53} -3 c_{46} -c_{64} -c_{35} +4 a_{44} +4 a_{55} \right) e^{145} \notag \\
        &+ \tfrac{1}{7} \left(2 a_{56} +2 c_{54} +2 a_{34} -2 c_{36} +5 a_{65} -5 c_{63} +5 a_{43} +5 c_{45} +2 b_{64} +2 b_{35} -2 b_{53} -2 b_{46} \right) e^{146} \notag \\
        &+ \tfrac{1}{4} \left(b_{65} +b_{43} -b_{56} -b_{34} -3 a_{64} +3 a_{53} -a_{46} +a_{35} -4 c_{44} -4 c_{33} \right) e^{156} \notag \\
        &+ \tfrac{1}{4} \left(-c_{56} +c_{43} +c_{65} -c_{34} +a_{63} +a_{54} +3 a_{45} +3 a_{36} -4 b_{33} -4 b_{44} \right) e^{234} \notag \\
        &+ \tfrac{1}{4} \left(b_{63} -b_{45} -3 b_{54} +3 b_{36} +c_{53} +c_{46} -c_{64} -c_{35} +4 a_{33} +4 a_{55} \right) e^{235} \notag \\
        &- \tfrac{1}{7} \left(-2 a_{56} -2 c_{54} +5 a_{34} +2 c_{36} -5 a_{65} -2 c_{63} +2 a_{43} +2 c_{45} +5 b_{64} +5 b_{35} +2 b_{53} +2 b_{46} \right) e^{236} \notag \\
        &+ \tfrac{1}{7} \left(-5 a_{56} +2 c_{54} +2 a_{34} -2 c_{36} -2 a_{65} +2 c_{63} +5 a_{43} -2 c_{45} +2 b_{64} +2 b_{35} +5 b_{53} +5 b_{46} \right) e^{245} \notag \\
        &+ \tfrac{1}{4} \left(3 b_{63} -3 b_{45} -b_{54} +b_{36} -c_{53} -c_{46} +c_{64} +c_{35} +4 a_{55} +4 a_{33} \right) e^{246} \notag \\
        &- \tfrac{1}{4} \left(c_{56} -c_{43} -c_{65} +c_{34} +3 a_{63} +3 a_{54} +a_{45} +a_{36} -4 b_{44} -4 b_{33} \right) e^{256} \notag \\
        &- \tfrac{1}{7} \left(2 a_{56} +2 c_{54} +2 a_{34} +5 c_{36} -2 a_{65} +2 c_{63} -2 a_{43} +5 c_{45} +2 b_{64} -5 b_{35} -2 b_{53} +5 b_{46} \right) e^{347} \notag \\
        &- \tfrac{1}{4} \left(b_{65} +b_{43} +3 b_{56} +3 b_{34} +a_{64} -a_{53} -a_{46} +a_{35} +4 c_{33} +4 c_{55} \right) e^{357} \notag \\
        &- \tfrac{1}{4} \left(c_{56} -c_{43} +3 c_{65} -3 c_{34} -a_{63} -a_{54} +a_{45} +a_{36} -4 b_{55} -4 b_{44} \right) e^{367} \notag \\
        &- \tfrac{1}{4} \left(-3 c_{56} +3 c_{43} -c_{65} +c_{34} -a_{63} -a_{54} +a_{45} +a_{36} +4 b_{44} +4 b_{55} \right) e^{457} \notag \\
        &+ \tfrac{1}{4} \left(-3 b_{65} -3 b_{43} -b_{56} -b_{34} +a_{64} -a_{53} -a_{46} +a_{35} -4 c_{55} -4 c_{33} \right) e^{467} \notag \\
        &- \tfrac{1}{7} \left(2 a_{56} -5 c_{54} +2 a_{34} -2 c_{36} -2 a_{65} -5 c_{63} -2 a_{43} -2 c_{45} -5 b_{64} +2 b_{35} +5 b_{53} -2 b_{46} \right) e^{567}. \notag
\end{align}
\end{proposition}
 
\indent The next result is a straightforward application of \eqref{omegas 1} of Lemma \ref{lemma: omegas}, where we put
\begin{align*}  
    k_1 &:= -  \tfrac{1}{12} ( a_{36} - a_{63} + a_{45} - a_{54} + c_{56} - c_{65} + c_{34} - c_{43} ), \\ 
        k_2 &:= -  \tfrac{1}{12} ( a_{64} - a_{46} + a_{35} - a_{53} + b_{43} - b_{34} + b_{65} - b_{56} ), \\
        k_7 &:= -  \tfrac{1}{12} (b_{63} - b_{36} + b_{54} - b_{45} + c_{46} - c_{64} + c_{53} - c_{35}).
\end{align*}  
\noindent so that $\tau_1 = k_1 e^1 + k_2 e^2 + k_7 e^7$.

\begin{corollary} \label{cor: iota N2}
    For each set of pairwise-commuting matrices $A$, $B$, $C \in \mathfrak{sl}(4, \mathbb{R})$, the following formula holds in $(\mathfrak{g}_{A,B,C}, \varphi)$: 
\begin{equation*}
    \iota_{\tau_1}(\varphi) = k_1 ( e^{27} + \omega_1 ) + k_2 ( - e^{17} + \omega_2) + k_7 ( e^{12} + \omega_7 ).
\end{equation*}    
\end{corollary} 

\begin{remark} \label{lo mas significativo}
    We note for future reference that Proposition \ref{prop: formulas N2} and Corollary \ref{cor: iota N2} imply that
\begin{equation*}
    \iota_{\tau_1}(\varphi), \tau_2 \in \Span \{e^{12}, e^{17}, e^{27}, e^{34}, e^{35}, e^{36}, e^{45}, e^{46}, e^{56}\}, 
\end{equation*} 
    \noindent as well as 
\begin{equation*}
    \tau_3 \in \Span\{e^{127}, e^{134}, e^{135}, e^{136}, e^{145}, e^{146}, e^{156}, e^{234}, e^{235}, e^{236}, e^{245}, e^{246}, e^{256}, e^{347}, e^{357}, e^{367}, e^{457}, e^{467}, e^{567}\}.
\end{equation*}      
\end{remark}

\subsection{Auxiliary formulas} 
\indent The following technical results aim to provide a better handling of the symmetric $2$-tensor $\tau_{27}$, which is a crucial step in the proof of Theorem \ref{thm: divergencia} below and ultimately allows us to establish two of the main results of this paper (Theorems \ref{thm: diag_divergenceless} and \ref{thm: adiag_divergenceless}). Lemma \ref{lemma: cuentoso} can be obtained by straightforward computations (see also Appendix \ref{section: maple 2023}), from which Corollary \ref{cor: cuentoso}. We assume only that $A$, $B$, $C \in \mathfrak{sl}(4, \mathbb{R})$ are pairwise-commuting and that $\varphi$ is given as in equation \eqref{eq: G2-estructura}.

\begin{lemma} \label{lemma: cuentoso}
The following formulas hold in $(\mathfrak{g}_{A,B,C}, \varphi)$: 
\begin{equation*}
    \iota_{e_j}(\varphi) =
    \begin{cases}
        + e^{27} + e^{35} - e^{46}, & j = 1. \\
        - e^{17} - e^{36} - e^{45}, & j = 2. \\
        + e^{47} - e^{15} + e^{26}, & j = 3. \\
        - e^{37} + e^{16} + e^{25}, & j = 4. \\
        + e^{67} + e^{13} - e^{24}, & j = 5. \\
        - e^{57} - e^{14} - e^{23}, & j = 6. \\
        + e^{12} + e^{34} + e^{56}, & j = 7. \\
    \end{cases}
\end{equation*}
\begin{equation*}
    \iota_{e_1}(\varphi) \wedge \iota_{e_j}(\varphi) = 
    \begin{cases}
        +e^{1257}-e^{3457}+e^{2356}-e^{1456}, & j = 3. \\
        -e^{1267}+e^{1356}+e^{3467}+e^{2456}, & j = 4. \\
        -e^{1237}+e^{3567}+e^{2345}-e^{1346}, & j = 5. \\
        +e^{1247}+e^{1345}-e^{4567}+e^{2346}, & j = 6.
    \end{cases}
\end{equation*}    
\begin{equation*}
    \iota_{e_2}(\varphi) \wedge \iota_{e_j}(\varphi) = 
    \begin{cases}
        -e^{1267}+e^{3467}-e^{1356}-e^{2456}, & j = 3. \\
        -e^{1257}+e^{2356}+e^{3457}-e^{1456}, & j = 4. \\
        +e^{1247}-e^{2346}-e^{4567}-e^{1345}, & j = 5. \\
        +e^{1237}-e^{3567}-e^{1346}+e^{2345}, & j = 6.
    \end{cases}
\end{equation*}    
\begin{equation*}
    \iota_{e_7}(\varphi) \wedge \iota_{e_j}(\varphi) = 
    \begin{cases}
        +e^{1247}-e^{1345}+e^{2346}+e^{4567}, & j = 3. \\
        -e^{1237}+e^{1346}+e^{2345}-e^{3567}, & j = 4. \\
        +e^{1267}+e^{3467}+e^{1356}-e^{2456}, & j = 5. \\
        -e^{1257}-e^{3457}-e^{1456}-e^{2356}, & j = 6.
    \end{cases}
\end{equation*}
\begin{equation*}
     \tfrac{1}{2} \iota_{e_j}(\varphi) \wedge \iota_{e_j}(\varphi)
    \begin{cases}
        +e^{1457} -e^{2467} +e^{1256}, & j = 3. \\
        +e^{1367} +e^{2357} +e^{1256}, & j = 4. \\
        +e^{1367} -e^{2467} +e^{1234}, & j = 5. \\
        +e^{1457} +e^{2357} +e^{1234}, & j = 6. \\
    \end{cases}
\end{equation*} 
\begin{equation*}
    \iota_{e_j}(\varphi) \wedge \iota_{e_{(9-j)}}(\varphi) = 
    \begin{cases}
        - e^{2347} + e^{1246} + e^{2567} + e^{1235} , & j = 3,6. \\
        - e^{2347} - e^{1246} + e^{2567} - e^{1235}, & j = 4,5. \\
    \end{cases}
\end{equation*} 
\end{lemma}

\begin{corollary} \label{cor: cuentoso}
    Define 
\begin{align*}
    \mathfrak{S}_1 &:= \Span\{e^{127}, e^{134}, e^{135}, e^{136}, e^{145}, e^{146}, e^{156}, e^{234}, e^{235}, e^{236}, e^{245}, e^{246}, e^{256}, e^{347}, e^{357}, e^{367}, e^{457}, e^{467}, e^{567}\}, \notag \\
    \mathfrak{S}_2 &:= \Span\{e^{134}, e^{136}, e^{145}, e^{156}, e^{234}, e^{235}, e^{246}, e^{256}, e^{357}, e^{367}, e^{457}, e^{467}\}, \notag \\
    \mathfrak{S}_3 &:= \Span\{e^{127}, e^{135}, e^{136}, e^{145},  e^{146}, e^{234}, e^{235}, e^{236}, e^{245}, e^{246}, e^{256}, e^{347}, e^{367}, e^{457}, e^{567}\}. \notag
\end{align*}    
    \noindent The following formulas hold in $(\mathfrak{g}_{A,B,C}, \varphi)$ for all $\chi_1 \in \mathfrak{S}_1$, $\chi_2 \in \mathfrak{S}_2$, and $\chi_3 \in \mathfrak{S}_3$:  
\begin{equation} \label{corolario cuentoso 1}
    \iota_{e_k}(\varphi) \wedge \iota_{e_i}(\varphi) \wedge \chi_1 = 0, \quad k = 1, 2, 7, \quad 3 \leq i \leq 6,
\end{equation}
\begin{equation} \label{corolario cuentoso 2}
    \iota_{e_n}(\varphi) \wedge \iota_{e_n}(\varphi) \wedge \chi_2 = 0, \quad 3 \leq n \leq 6, 
\end{equation}  
\begin{equation} \label{corolario cuentoso 3}
    \iota_{e_m}(\varphi) \wedge \iota_{e_{(9-m)}}(\varphi) \wedge \chi_3 = 0, \quad m \leq l \leq 6.
\end{equation} 
\end{corollary}

\indent Following the definition of $\tau_{27}$ given in equation \eqref{tau_27}, the next result is an immediate consequence of Remark \ref{lo mas significativo} and equation \eqref{corolario cuentoso 1}.

\begin{corollary} \label{cor: tau27 en el caso general}
    The following formulas hold in $(\mathfrak{g}_{A,B,C}, \varphi)$:
\begin{equation*}
    \tau_{27}(e_k, e_i) = 0 \quad \text{for all $k = 1, 2, 7$ and $3 \leq i \leq 6$}.
\end{equation*}    
\end{corollary}  

\subsection{Formulas for the Levi-Civita connection and the Ricci operator}
\indent We now turn to the goal of finding the Levi-Civita connection of the underlying Riemannian structure in $(\mathfrak{g}_{A,B,C}, \varphi)$. Special notation is to be introduced in order to state the result more briefly:
\begin{itemize}
    \item Given $Z \in \mathfrak{a}$ of the form $Z = \lambda_1 \, e_1 + \lambda_2 \, e_2 + \lambda_7 \, e_7$ for some $\lambda_1$, $\lambda_2$, $\lambda_7 \in \mathbb{R}$, we define
\begin{equation*}
    M_Z := \lambda_1 \, B + \lambda_2 \, C + \lambda_7 \, A.
\end{equation*}   
\noindent Note that if $Z \in \mathfrak{a}$ and $W \in \mathfrak{n}$ then $[Z,W] = M_Z W$. 
    \item We call $S(M)$ and $A(M)$ the symmetric and antisymmetric parts of a matrix $M \in \mathfrak{gl}(4, \mathbb{R})$ respectively; that is,
\begin{equation*}
    S(M) :=  \tfrac{M + M^{\intercal}}{2}, \quad A(M) :=  \tfrac{M - M^{\intercal}}{2}.
\end{equation*} 
    \item We define 
\begin{equation} \label{eq: increible pero necesito un label}
    D^l = B \delta_{l1} + C \delta_{l2} + A \delta_{l7},
\end{equation}    
    \noindent where $\delta_{ij}$ is the Kronecker delta and $A$, $B$, $C \in \mathfrak{sl}(4, \mathbb{R})$ are the matrices appearing in the definition of $\mathfrak{g}_{A,B,C}$. Note that we are \textbf{not} performing a summation over the index $l$. 
\end{itemize}

\begin{proposition} \label{prop: LC in g_ABC}
    Let $A$, $B$, $C \in \mathfrak{sl}(4, \mathbb{R})$ be pairwise-commuting. The Levi-Civita connection $\nabla$ of the underlying Riemannian structure in $(\mathfrak{g}_{A,B,C}, \varphi)$ is given by
\begin{equation*}
    \nabla_X Y =
    \left\{
        \begin{array}{cl}
            0 & X \in \mathfrak{a}, \, Y \in \mathfrak{a}.      \\
            \phantom{+} A(M_X) Y & X \in \mathfrak{a}, \, Y \in \mathfrak{n}.     \\
            - S(M_Y) X & X \in \mathfrak{n}, \, Y \in \mathfrak{a}.     \\
            \sum_{l=1,2,7} \langle S(D^l) X,Y \rangle \, e_l & X \in \mathfrak{n}, \, Y \in \mathfrak{n}.  \\
        \end{array} 
    \right.
\end{equation*}
\end{proposition}

\begin{proof}
    We consider four separate cases.  
    
    Let $X$, $Y \in \mathfrak{a}$. Then $[X,Y] = 0$, as $\mathfrak{a}$ is an abelian subalgebra. Therefore, equation \eqref{eq: U(X,Y)} implies that $U(X,Y)$ has no components in $\mathfrak{a}$. In addition, $U(X,Y)$ has no components in $\mathfrak{n}$, for if $Z \in \mathfrak{n}$ then $[Z,X] \in \mathfrak{n}$ and $[Z,Y] \in \mathfrak{n}$ because $\mathfrak{n}$ is an ideal, and then
\begin{equation*}
    2 \, \langle U(X,Y),Z \rangle = \langle [Z,X],Y \rangle + \langle [Z,Y],X \rangle = 0,
\end{equation*}    
    \noindent where the last equality holds because $\mathfrak{a}$ and $\mathfrak{n}$ are orthogonal complements in $\mathfrak{g}_{A,B,C}$. Thus, equation \eqref{eq: nabla linealizado} implies that $\nabla_X Y = 0$ if $X$, $Y \in \mathfrak{a}$. 
    
    Let $X \in \mathfrak{a}$ and $Y \in \mathfrak{n}$. As equation \eqref{eq: nabla linealizado} implies that $\nabla_X Y = \frac{M_X Y}{2} + U(X,Y)$, it suffices to show that $U(X,Y) = - \frac{M_X^{\intercal} Y}{2}$. To such end, let $W \in \mathfrak{a}$ and $Z \in \mathfrak{n}$ be arbitrary. Note that $[W,X] = 0$ and $[Z,Y] = 0$ because both $\mathfrak{a}$ and $\mathfrak{n}$ are abelian. In addition, as $[W,Y] \in \mathfrak{n}$ for $\mathfrak{n}$ is an ideal, it follows that $\langle [W,Y],X \rangle = 0$ due to $\mathfrak{a}$ and $\mathfrak{n}$ being orthogonal complements in $\mathfrak{g}_{A,B,C}$. Thus, equation \eqref{eq: U(X,Y)} implies that
    \begin{equation*}
        2 \, \langle U(X,Y), W \rangle = \langle [W,X],Y \rangle + \langle [W,Y],X \rangle = 0;
    \end{equation*}
    \noindent this establishes that $U(X,Y) \in \mathfrak{n}$. Furthermore, as $[Z,X] = - M_X Z$, the same equations show that
    \begin{equation*}
        2 \, \langle U(X,Y), Z \rangle = \langle [Z,X],Y \rangle + \langle [Z,Y],X \rangle = - \langle M_X Z,Y \rangle = - \, \langle M_X^\intercal Y,Z \rangle.
    \end{equation*}
    \noindent We can write the last equality $\langle U(X,Y), Z \rangle = - \frac{\langle M_X^{\intercal} Y, Z\rangle}{2}$. Combining the last two results, we get $U(X,Y) = - \frac{M_X^{\intercal} Y}{2}$. 
    
    Let $X \in \mathfrak{n}$ and $Y \in \mathfrak{a}$. From the torsion-free property of $\nabla$ and the previous paragraph, we see that
    \begin{equation*}
        \nabla_X Y = \nabla_Y X + [X,Y] = A(M_Y) X - M_Y X = - S(M_Y) X.  
    \end{equation*}
    
    Let $X,Y \in \mathfrak{n}$. Then $[X,Y] = 0$, as $\mathfrak{n}$ is abelian. Therefore, equation \eqref{eq: U(X,Y)} implies that $U(X,Y)$ has no components in $\mathfrak{n}$. Let $Z \in \mathfrak{a}$, e.g. $Z = e_7$. We see from equation \eqref{eq: U(X,Y)} that
    \begin{equation*}
        \begin{split}
        2 \, \langle U(X,Y),e_7 \rangle &= \langle [e_7,X],Y \rangle + \langle [e_7,Y],X \rangle = \langle AX,Y \rangle + \langle AY,X \rangle = \langle (A + A^\intercal)X,Y \rangle; 
        \end{split}
    \end{equation*}
    \noindent i.e., that $\langle U(X,Y), e_7 \rangle = \langle S(A)X, Y \rangle$. Doing the same for $Z = e_1$ y $Z = e_2$, we find that
\begin{equation*}
    U(X,Y) = \langle S(B) X,Y \rangle \, e_1 + \langle S(C) X,Y \rangle \, e_2 + \langle S(A) X,Y \rangle \, e_7. \qedhere
\end{equation*}    
\end{proof}
\begin{remark}
    A similar result to Proposition \ref{prop: LC in g_ABC} is valid in any Lie algebra $\mathfrak{g}$ that splits as an orthogonal direct sum $\mathfrak{a} \oplus \mathfrak{n}$ with respect to an inner product $\langle \cdot, \cdot \rangle$, where $\mathfrak{a}$ is an abelian subalgebra and $\mathfrak{n}$ is an abelian ideal. 
\end{remark}

\noindent We recall the following result for further reference. 

\begin{proposition} \cite[Section 2.5]{N} \label{Ric in g_ABC}
    Let $A$, $B$, $C \in \mathfrak{sl}(4, \mathbb{R})$ be pairwise-commuting. The Ricci operator $\mathrm{Ric}:\mathfrak{g}_{A,B,C} \to \mathfrak{g}_{A,B,C}$ of the underlying Riemannian structure in $(\mathfrak{g}_{A,B,C}, \varphi)$ is given by
\begin{align*}
    \mathrm{Ric} \vert_{\mathfrak{a} \times \mathfrak{n}} &= 0, \\
    \mathrm{Ric} \vert_{\mathfrak{n} \times \mathfrak{n}} &= \tfrac{1}{2} ([A, A^{\intercal}] + [B, B^{\intercal}] + [C, C^{\intercal}]), \\
    \mathrm{Ric} \vert_{\mathfrak{a} \times \mathfrak{a}} &= -
    \begin{pmatrix}
        \tr(S(A)^2) & \tr(S(A) B) & \tr(S(A) C) \\
        \tr(S(A) B) & \tr(S(B)^2) & \tr(S(B) C) \\
        \tr(S(A) C) & \tr(S(B) C) & \tr(S(C)^2)
    \end{pmatrix}.
\end{align*}    
\end{proposition}

\subsection{Formula for the divergence of the torsion}

\indent Remark \ref{lo mas significativo}, Corollary \ref{cor: tau27 en el caso general}, and Proposition \ref{prop: LC in g_ABC} to give a simple expression for the divergence of the full torsion tensor of $(\mathfrak{g}_{A,B,C}, \varphi)$, given by equation \eqref{eq: div alg lie}.  

\begin{theorem} \label{thm: divergencia}
    Let $A$, $B$, $C \in \mathfrak{sl}(4, \mathbb{R})$ be pairwise-commuting. The divergence $\DV T$ of the full torsion tensor $T$ of $(\mathfrak{g}_{A,B,C}, \varphi)$ satisfies
\begin{equation*}
    \langle \DV T, e_j \rangle = 
    \begin{cases}
        -\sum_{3 \leq n \leq 6} (D^j)_{nn} \tau_{27}(e_n, e_n) + \sum_{3 \leq i \neq l \leq 6} S(D^j)_{il} \tau_{27}(e_i, e_l) & j = 1,2,7,  \\
        \phantom{\sum_{3 \leq n \leq 6} (D^j)_{nn} \tau_{27}(e_n, e_n)} \quad \; 0 & 3 \leq j \leq 6,
    \end{cases}
\end{equation*}
    \noindent with $D^j$ as in equation \eqref{eq: increible pero necesito un label}. 
\end{theorem} 
\begin{proof}
     Recall from equation \eqref{eq: div alg lie} that
\begin{equation*} 
    \begin{split}
        \langle \DV T, e_j \rangle &= - \sum_{1 \leq i \leq 7} T(\nabla_{e_i} e_i, e_j) - \sum_{1 \leq i \leq 7} T(e_i, \nabla_{e_i} e_j).
    \end{split}
\end{equation*}
\noindent We claim that Proposition \ref{prop: LC in g_ABC} and $A$, $B$, $C$ being traceless ensure that the first term is identically zero. Indeed, for $1 \leq j \leq 7$, we have
\begin{equation*}
    \begin{split}
        \sum_{1 \leq i \leq 7} T(\nabla_{e_i} e_i, e_j) &= \sum_{3 \leq i \leq 6} T(\nabla_{e_i} e_i, e_j) = \sum_{3 \leq i \leq 6} T \left( \sum_{k=1,2,7} \langle S(D^k) e_i, e_i\rangle e_k, e_j \right)\\
        &= \sum_{3 \leq i \leq 6} \sum_{k=1,2,7} S(D^k)_{ii} T(e_k, e_j) = \sum_{k=1,2,7} \tr(S(D^k)) T(e_k, e_j)\\
        &= \sum_{k=1,2,7} \tr(D^k) T(e_k, e_j) = \sum_{k=1,2,7} 0 \cdot T(e_k, e_j) = 0.
    \end{split}
\end{equation*}
    \noindent On the other hand, Proposition \ref{prop: LC in g_ABC} ensures that for $j = 1, 2, 7$ we have
\begin{equation*}
    \begin{split}
    \sum_{1 \leq i \leq 7} T(e_i, \nabla_{e_i} e_j) &= \sum_{3 \leq i \leq 6} T(e_i, \nabla_{e_i}  e_j) = - \sum_{3 \leq i \leq 6} T(e_i, S(D^j) e_i) = - \sum_{3 \leq i,l \leq 6}  S(D^j)_{il} T(e_i, e_l) \\
    &= \sum_{3 \leq n \leq 6} (D^j)_{nn} T(e_n, e_n) + \sum_{3 \leq i \neq l \leq 6} S(D^j)_{il} T(e_i, e_l),
    \end{split}
\end{equation*} 
    \noindent whereas for $3 \leq j \leq 6$ we have
\begin{equation*}
    \begin{split}
    \sum_{1 \leq i \leq 7} T(e_i, \nabla_{e_i} e_j) &= \sum_{i=1,2,7} T(e_i, A(D^i) e_j) + \sum_{3 \leq i \leq 6} T\left(e_i, \sum_{k=1,2,7} \langle S(D^k) e_i, e_j \rangle \, e_k \right) \\
    &= \sum_{i=1,2,7} T(e_i, A(D^i) e_j) + \sum_{3 \leq i \leq 6} \sum_{k=1,2,7} S(D^k)_{ij} T(e_i,e_k) \\
    &= \sum_{k=1,2,7} \left( T(e_k, A(D^k) e_j) + \sum_{3 \leq i \leq 6} S(D^k)_{ij} T(e_i, e_k) \right).  
    \end{split}
\end{equation*}
    \noindent Let us recall now equation \eqref{eq: full torsion tensor form}, which states that 
\begin{equation} \label{eq: full torsion tensor form in g}
    T(u,v) =  \tfrac{1}{4} \tau_0 g(u,v) - \iota_{\tau_1}(\varphi)(u,v) -  \tfrac{1}{2} \tau_2(u,v) - \tau_{27}(u,v), \quad u,v \in \mathfrak{g}_{A,B,C},
\end{equation}    
    \noindent and equation \eqref{tau_27} defining the symmetric $2$-tensor $\tau_{27}$. When $u = v = e_n$ for $3 \leq n \leq 6$ we find that the first term in equation \eqref{eq: full torsion tensor form in g} is $ \tfrac{\tau_0}{4}$ because $\{e_j \, \vert \, 1 \leq j \leq 7 \}$ is chosen to be an orthonormal basis, and that the second and third terms are zero because they are skew-symmetric, but the last term is non-zero in general; moreover, given that $\sum_{3 \leq n \leq 6} (D^j)_{nn}  \tfrac{\tau_0}{4} = 0$ because of tracelessness of $D = A$, $B$, $C$, we have that
\begin{equation*}
    \sum_{3 \leq n \leq 6} (D^j)_{nn} T(e_n, e_n) = - \sum_{3 \leq n \leq 6} (D^j)_{nn} \tau_{27} (e_n, e_n).
\end{equation*}    
    \noindent Similarly, when $u = e_i$ and $v = e_l$ for $3 \leq i \neq l \leq 6$, we find that the first term in equation \eqref{eq: full torsion tensor form in g} is zero because $\{e_j \, \vert \, 1 \leq j \leq 7 \}$ is chosen to be an orthonormal basis; moreover, the fact that $S(D^j)_{il}$ is symmetric in the indices $i$ and $l$ while the second and third terms in equation \eqref{eq: full torsion tensor form in g} are skew-symmetric entails that
\begin{equation*}
    \sum_{3 \leq i \neq l \leq 6} S(D^j)_{il} T(e_i, e_l) = - \sum_{3 \leq i \neq l \leq 6} S(D^j)_{il} \tau_{27}(e_i, e_l).
\end{equation*}      
    \noindent We can then summarize the findings of the last paragraph as
\begin{equation*}
    \langle \DV T, e_j \rangle = -\sum_{3 \leq n \leq 6} (D^j)_{nn} \tau_{27}(e_n, e_n) + \sum_{3 \leq i \neq l \leq 6} S(D^j)_{il} \tau_{27}(e_i, e_l) \quad \text{for $j = 1, 2, 7$}.
\end{equation*}    
    \noindent On the other hand, when $u = e_i$ and $v = e_k$ for $3 \leq i \leq 6$ and $k = 1,2,7$ we find that the first term in equation \eqref{eq: full torsion tensor form in g} is zero because $\{e_j \, \vert \, 1 \leq j \leq 7 \}$ is chosen to be an orthonormal basis, the second and third terms are zero because of Remark \ref{lo mas significativo}, and the last term is zero because of Remark \ref{lo mas significativo} and Corollary \ref{cor: tau27 en el caso general}. This means that $T(e_i, e_k) = 0$ for $3 \leq i \leq 6$ and $k = 1,2,7$; moreover, it also implies that $T(e_k, A(D^k) e_j) = 0$ for $k = 1, 2, 7$ and $3 \leq j \leq 6$, for $A(D^k) e_j \in \mathfrak{n} = \Span \{e_3, \pts, e_6\}$ and $T$ is bilinear. Therefore, for $3 \leq j \leq 6$ we have 
\begin{equation*}
    \sum_{1 \leq i \leq 7} T(e_i, \nabla_{e_i} e_j) = \sum_{k=1,2,7} \left( T(e_k, A(D^k) e_j) + \sum_{3 \leq i \leq 6} S(D^k)_{ij} T(e_i, e_k) \right) = 0, 
\end{equation*}    
    \noindent and thus,
\begin{equation*}
    \langle \DV T, e_j \rangle = 0 \quad \text{for $3 \leq j \leq 6$}. \qedhere 
\end{equation*}    
\end{proof} 

\indent Theorem \ref{thm: divergencia} reveals that in order to have an explicit formula for $\DV T$ we only need to compute 
\begin{equation} \label{eq: los tau27 clave}
    \tau_{27}(e_n, e_n) \quad \text{for $3 \leq n \leq 6$}, \quad \tau_{27}(e_i, e_l) \quad \text{for $3 \leq i \neq j \leq 6$}, 
\end{equation}
\noindent a task rendered possible by equation \eqref{tau3 N2}. It appears that little gain is made in performing such amount of computations in this general context, so we omit it. It remains unclear how to characterize when the full torsion tensor $T$ is divergence-free in terms of conditions on the coefficients of $A$, $B$, and $C$, and even more so which of those conditions lead to non-equivalent $\mathrm{G}_2$-structures. Going forward, we will restrict ourselves to the three broad cases that proved tractable: The case in which $A$, $B$, $C$ are all skew-symmetric, the case in which $A$, $B$, $C$ are all diagonal (which turns out to be equivalent to the symmetric case, see Proposition \ref{prop: A,B,C sim compt}), and the case in which $A$, $B$, $C$ are all antidiagonal. Recall that a matrix $D = [d_{ij}]$ for $3 \leq i,j \leq 6$ is called \textit{antidiagonal} if it is of the form
\begin{equation} \label{eq: definicion antidiagonal}
   D = \begin{psmallmatrix}
     0 & 0 & 0 & d_{36} \\
     0 & 0 & d_{45} & 0 \\
     0 & d_{54} & 0 & 0 \\
     d_{63} & 0 & 0 & 0 
   \end{psmallmatrix}.
\end{equation}
\noindent Generally speaking, these cases are indeed non-equivalent as a consequence of Proposition \ref{Ric in g_ABC}, for we have the following:
\begin{itemize}
    \item In the skew symmetric case, $\mathrm{Ric}$ is identically zero (in fact, $\nabla$ is flat in this case; see Proposition \ref{antisimetrico es plano}).
    \item In the diagonal (or symmetric) case, both $\mathrm{Ric} \vert_{\mathfrak{a} \times \mathfrak{n}}$ and $\mathrm{Ric} \vert_{\mathfrak{n} \times \mathfrak{n}}$ are identically zero, whereas $\mathrm{Ric} \vert_{\mathfrak{a} \times \mathfrak{a}}$ is generally not.
    \item In the antidiagonal case, both $\mathrm{Ric} \vert_{\mathfrak{a} \times \mathfrak{n}}$ and $\mathrm{Ric} \vert_{\mathfrak{a} \times \mathfrak{a}}$ are identically zero, whereas $\mathrm{Ric} \vert_{\mathfrak{n} \times \mathfrak{n}}$ is generally not at least when $A$, $B$, $C$ are not all skew-symmetric or diagonal/symmetric.
\end{itemize} 

\section{Skew-symmetric case} \label{skew case}

\indent This case is of particular interest due in part to the following well-known result.

\begin{proposition} \cite[Theorem 1.5]{Mi} \label{antisimetrico es plano}
    If $A$, $B$, $C \in \mathfrak{sl}(4, \mathbb{R})$ are skew-symmetric pairwise-commuting then the metric $\langle \cdot, \cdot \rangle$ induced in $\mathfrak{g}_{A,B,C}$ by $\varphi$ is flat.
\end{proposition}

\begin{remark}
    Proposition \ref{antisimetrico es plano} does \textbf{not} mean that $\nabla \varphi = 0$, only that the curvature tensor of the metric induced by $\varphi$ vanishes. See Proposition \ref{prop: taus skew} below. 
\end{remark} 

Flatness (or, more generally, Ricci-flatness) is a relevant assumption (though not crucial) in the context of the isometric flow of $\mathrm{G}_2$-structures, for torsion-free $\mathrm{G}_2$-structures are known to be Ricci-flat (see \cite[Propositions 10.1.3 and 10.1.5]{J}) and the isometric flow preserves the underlying metric. Hence the importance of Proposition \ref{antisimetrico es plano}. 

\indent Another remark  worth mentioning is that three pairwise-commuting skew-symmetric matrices are always a linearly-dependent set, for the dimension of the Cartan subalgebra\footnote{By \textit{Cartan subalgebra} of $\mathfrak{so}(3)$ we mean a maximal abelian Lie subalgebra consisting of diagonalizable elements of the complexification of $\mathfrak{so}(3)$; it is well known that any two Cartan subalgebras are isomorphic.} of $\mathfrak{so}(3)$ is 2, and so there are no compatible triples $A$, $B$, $C \in \mathfrak{sl}(4, \mathbb{R})$ that are skew-symmetric.

\indent According to Theorem \ref{thm: divergencia}, if all three matrices $A$, $B$, $C \in \mathfrak{sl}(4, \mathbb{R})$ are skew-symmetric then the divergence of the full torsion tensor $T$ of $(\mathfrak{g}_{A,B,C}, \varphi)$ is identically zero, and trivially so as we have $(D^j)_{nn} = 0 = S(D^j)_{il}$ for all $j = 1, 2, 7$ and $3 \leq i, l \leq 6$. We state this result explicitly.

\begin{theorem} \label{antisim_divergenceless}
    Let $A$, $B$, $C \in \mathfrak{sl}(4, \mathbb{R})$ be skew-symmetric. Then $(\mathfrak{g}_{A,B,C}, \varphi)$ is a $\mathrm{G}_2$-structure with divergence-free full torsion tensor. 
\end{theorem}
\begin{remark}
    Theorem \ref{antisim_divergenceless} is consistent with \cite[Theorem 4.8]{T}, in which a much larger class of $7$-dimensional solvable Lie groups is found to admit a $\mathrm{G}_2$-structure with divergence-free full torsion tensor that is compatible with a given flat left-invariant metric.
\end{remark}
 
\indent We adapt Proposition \ref{prop: formulas N2} to the diagonal case for the sake of completeness; see Appendix \ref{section: maple 2023}. Incidentally, this illustrates that non-torsion-free $\mathrm{G}_2$-structures may induce flat riemannian metrics. 

\begin{proposition} \label{prop: taus skew}
    The torsion forms $\tau_0$, $\tau_1$, $\tau_2$, and $\tau_3$ of $(\mathfrak{g}_{A,B,C}, \varphi)$ in the skew-symmetric case are given by the following formulas: 
\begin{align*}
    \tau_0 =& \phantom{+} \; \, \tfrac{4}{7} (a_{34} + a_{56} + b_{35} + b_{64} + c_{54} + c_{63}), \\
    \tau_1 =& -  \tfrac{1}{6} ( a_{36} + a_{45} + c_{56} + c_{34}) e^1 -  \tfrac{1}{6} (a_{64} + a_{35} +  b_{43} + b_{65} ) e^2 -  \tfrac{1}{6} (b_{63} + b_{54}+ c_{46} + c_{53}) e^7, \\
    \tau_2 =& \phantom{+} \; \, \tfrac{2}{3} (b_{45} + b_{36} + c_{35} + c_{64}) e^{12} +  \tfrac{2}{3} (a_{64} + a_{35} + b_{65} + b_{43} ) e^{17} +  \tfrac{2}{3} ( a_{54} + a_{63} + c_{65} + c_{43}) e^{27} \label{tau2 antisim} \\
    &+  \tfrac{1}{3} (c_{46} - c_{35} - b_{45} - b_{36}) e^{34} +  \tfrac{1}{3} (- a_{54} + a_{36} + c_{56} + c_{34}) e^{35} +  \tfrac{1}{3} (- a_{64} - a_{35} - b_{65} + b_{34}) e^{36} \notag \\
    &+  \tfrac{1}{3} (a_{46} + a_{53} + b_{56} - b_{43}) e^{45} +  \tfrac{1}{3} (a_{63} - a_{45} + c_{65} + c_{43}) e^{46} +  \tfrac{1}{3} (c_{53} - c_{64} + b_{63} + b_{54}) e^{56}, \notag \\ 
    \tau_3 =& \phantom{+} \; \, \tfrac{4}{7} (a_{65}-c_{63}+a_{43}-c_{54}+b_{53}-b_{64}) e^{127} + \tfrac{1}{2} (b_{65}+b_{43}-a_{64}+a_{53}) e^{134} \\
    &+ \tfrac{1}{7} (-3 a_{65}+3 c_{63}-3 a_{43}+3 c_{54}+4 b_{53}-4 b_{64}) e^{135}  + \tfrac{1}{2} (b_{63}+b_{54}-c_{53}+c_{64}) e^{136} \notag \\
    &+ \tfrac{1}{2} (b_{63}+b_{54}-c_{53}+c_{64}) e^{145} + \tfrac{1}{7} (3 a_{65}-3 c_{63}+3 a_{43}-3 c_{54}-4 b_{53}+4 b_{64}) e^{146} \notag \\
    &+ \tfrac{1}{2} (b_{65}+b_{43}-a_{64}+a_{53}) e^{156} + \tfrac{1}{2} (c_{65}+c_{43}-a_{54}-a_{63}) e^{234} \notag \\
    &+ \tfrac{1}{2} (-b_{63}-b_{54}+c_{53}-c_{64}) e^{235} + \tfrac{1}{7} (3 a_{65}+4 c_{63}+3 a_{43}+4 c_{54}+3 b_{53}-3 b_{64}) e^{236} \notag \\
    &+ \tfrac{1}{7} (3 a_{65}+4 c_{63}+3 a_{43}+4 c_{54}+3 b_{53}-3 b_{64}) e^{245} + \tfrac{1}{2} (b_{63}+b_{54}-c_{53}+c_{64}) e^{246} \notag \\
    &+ \tfrac{1}{2} (c_{65}+c_{43}-a_{54}-a_{63}) e^{256} + \tfrac{1}{7} (4 a_{65}+3 c_{63}+4 a_{43}+3 c_{54}-3 b_{53}+3 b_{64}) e^{347} \notag \\
    &+ \tfrac{1}{2} (b_{65}+b_{43}-a_{64}+a_{53}) e^{357} + \tfrac{1}{2} (-c_{65}-c_{43}+a_{54}+a_{63}) e^{367} \notag \\
    &+ \tfrac{1}{2} (-c_{65}-c_{43}+a_{54}+a_{63}) e^{457} + \tfrac{1}{2} (-b_{65}-b_{43}+a_{64}-a_{53}) e^{467} \notag \\
    &+ \tfrac{1}{7} (4 a_{65}+3 c_{63}+4 a_{43}+3 c_{54}-3 b_{53}+3 b_{64}) e^{567}. \notag
\end{align*}
\end{proposition} 

\section{Diagonal case} \label{diag case}

\indent This case is of particular interest due in part to the following result.
\begin{proposition} \cite[Lemma 4.5]{KL} \label{prop: A,B,C sim compt}
    If $A$, $B$, $C \in \mathfrak{sl}(4, \mathbb{R})$ are symmetric pairwise-commuting then there exist diagonal $A_0$, $B_0$, $C_0 \in \mathfrak{sl}(4, \mathbb{R})$ such that $(G_{A,B,C}, \varphi)$ is equivariantly equivalent to $(G_{A_0,B_0,C_0}, \varphi)$, where $\varphi$ is as given in equation \eqref{positiva}. 
\end{proposition}
\indent Proposition \ref{prop: A,B,C sim compt} shows that the diagonal case is in some sense the ``opposite"\! of the skew-symmetric case; moreover, it ensures that it is enough to consider only diagonal triples $A$, $B$, $C \in \mathfrak{sl}(4, \mathbb{R})$ to fully understand the case in which all three $A$, $B$, $C \in \mathfrak{sl}(4, \mathbb{R})$ are symmetric and pairwise-commuting. Note that any triple of diagonal matrices is trivially pairwise-commuting, so tracelessness is the only extra condition we need to impose.  

\indent According to Theorem \ref{thm: divergencia}, if all three matrices $A$, $B$, $C$ are diagonal then the divergence of the full torsion tensor $T$ of $(\mathfrak{g}_{A, B, C}, \varphi)$ is
\begin{equation} \label{eq: divergencia formula diag}
    \langle \DV T, e_j \rangle = 
    \begin{cases}
        - \sum_{3 \leq n \leq 6} (D^j)_{nn} \, \tau_{27}(e_n, e_n) & j = 1,2,7,  \\
        \phantom{- \sum_{3 \leq n \leq 6} (D^j)_{nn}} \; 0 & 3 \leq j \leq 6.
    \end{cases} 
\end{equation}
\noindent To proceed, we first adapt Proposition \ref{prop: formulas N2} to the diagonal case; see Appendix \ref{section: maple 2023}. 

\begin{proposition} \label{prop: taus diag}
    The torsion forms $\tau_0$, $\tau_1$, $\tau_2$, and $\tau_3$ of $(\mathfrak{g}_{A,B,C}, \varphi)$ in the diagonal case are given by the following formulas:
\begin{align}
    \tau_0 =& \; 0, \label{tau_0 diag} \\
    \tau_1 =& \; 0, \label{tau_1 diag} \\
    \tau_2 =& - (a_{33} + a_{44}) e^{34} - (b_{33} + b_{55}) e^{35} - ( c_{44} + c_{55}) e^{36} \label{tau_2 diag} \\
    &+ (c_{44} + c_{55}) e^{45} - (b_{33} - b_{55}) e^{46} + (a_{33} + a_{44}) e^{56}, \notag \\
    \tau_3 =& \phantom{+} \; \, (c_{33} + c_{44}) e^{134} - (a_{44} + a_{55}) e^{136} + (a_{44} + a_{55}) e^{145} \label{tau_3 diag}     \\ 
    &- (c_{33} + c_{44}) e^{156} - (b_{33} + b_{44}) e^{234} + (a_{33} + a_{55}) e^{235} \notag \\  
    &+ (a_{33} + a_{55}) e^{246} + (b_{33} + b_{44}) e^{256} - (c_{33} + c_{55}) e^{357} \notag \\ 
    &+ (b_{44} + b_{55}) e^{367} - (b_{44} + b_{55}) e^{457} - (c_{33} + c_{55}) e^{467}, \notag
\end{align}   
\end{proposition}

\indent Given that $\tau_3 \in \Span \{e^{134}, e^{136}, e^{145}, e^{156}, e^{234}, e^{235}, e^{246}, e^{256}, e^{357}, e^{367}, e^{457}, e^{467}\}$, equation \eqref{corolario cuentoso 2} from Corollary \eqref{cor: cuentoso} entails the following result. 
\begin{corollary} \label{cor: tau27 en el caso diagonal}
    The symmetric $2$-tensor $\tau_{27}$ of $(\mathfrak{g}_{A,B,C}, \varphi)$ in the diagonal case satisfies 
\begin{equation*}
    \tau_{27}(e_n, e_n) = 0 \quad \text{for all $3 \leq n \leq 6$}.
\end{equation*}    
\end{corollary}
\indent We apply Corollary \ref{cor: tau27 en el caso diagonal} to equation \eqref{eq: divergencia formula diag} to establish the main result of this Section.

\begin{theorem} \label{thm: diag_divergenceless}
    Let $A$, $B$, $C \in \mathfrak{sl}(4, \mathbb{R})$ be diagonal. Then $(\mathfrak{g}_{A,B,C}, \varphi)$ is a $\mathrm{G}_2$-structure with divergence-free full torsion tensor. 
\end{theorem} 

\indent Theorem \ref{thm: diag_divergenceless} can be enhanced as a result of Proposition \ref{prop: A,B,C sim compt}. 

\begin{corollary} \label{sim divergenceless}
    Let $A$, $B$, $C \in \mathfrak{sl}(4, \mathbb{R})$ be symmetric and pairwise-commuting. Then $(\mathfrak{g}_{A,B,C}, \varphi)$ is a $\mathrm{G}_2$-structure with divergence-free full torsion tensor.
\end{corollary}

\begin{remark} \label{obs: closed are divergence-free}
    Combining equations \eqref{tau_0 diag} and \eqref{tau_1 diag} with equation \eqref{eq: phi phis} yields $\tau_2 = - \star d\psi$ and $\tau_3 = \star d \varphi$, from which follows that $0 = \star (d \star \tau_2)$ and $0 = \star (d \star \tau_3)$. From this we get that $\tau_2$ and $\tau_{27}$ are divergence-free tensors (see \cite[Proposition 2.11]{G1}), and, by means of equation \eqref{eq: full torsion tensor form}, that $\DV T = 0$. Notice that Theorem \ref{thm: divergencia} was not used in this reasoning. The fact that a $\mathrm{G}_2$-structure with $\tau_0 = 0$ and $\tau_1 = 0$ has divergence-free full torsion tensor was observed previously in \cite[Theorem 4.3]{G2}. 
\end{remark}


\section{Antidiagonal case} \label{adiag case} 

\indent It is not known whether an analogous result to Proposition \ref{prop: A,B,C sim compt} of Section \ref{diag case} holds in the antidiagonal case. Note that in order to produce non-equivalent examples of $\mathrm{G}_2$-structures with divergence-free full torsion tensor to those found in the previous sections we need to consider only triples $A$, $B$, $C \in \mathfrak{sl}(4, \mathbb{R})$ such that not all three matrices are skew-symmetric or symmetric. We also note that it is not guaranteed that any election of a triple $A$, $B$, $C \in \mathfrak{sl}(4, \mathbb{R})$ of antidiagonal matrices pairwise commute. On the contrary, this imposes a series of restrictions on such matrices. Those restrictions, however, play no role in what follows, and therefore any discussion regarding them is omitted.

\indent According to Theorem \ref{thm: divergencia}, if all three matrices $A$, $B$, $C$ are antidiagonal then the divergence of the full torsion tensor $T$ of $(\mathfrak{g}_{A,B,C}, \varphi)$ is
\begin{equation} \label{eq: divergencia formula adiag}
    \langle \DV T, e_j \rangle = 
    \begin{cases}
        \sum_{3 \leq m \leq 6} S(D^j)_{m, 9-m} \, \tau_{27}(e_m, e_{9 - m}) & j = 1,2,7,  \\
        \phantom{\sum_{3 \leq m \leq 6} S(D^j)_{m, 9-m}} \; \, 0 & 3 \leq j \leq 6.
    \end{cases} 
\end{equation}
\noindent As in the previous Sections, we first adapt Proposition \ref{prop: formulas N2} to the antidiagonal case; see Appendix \ref{section: maple 2023}. 

\begin{proposition} \label{prop: taus adiag}
    The torsion forms $\tau_0$, $\tau_1$, $\tau_2$, and $\tau_3$ of $(\mathfrak{g}_{A,B,C}, \varphi)$ in the antidiagonal case are given by the following formulas:
\begin{align*}
    \tau_0 =& \phantom{+} \; \, \tfrac{2}{7} (c_{54} - c_{45} + c_{63} - c_{36}),  \\
    \tau_1 =& \phantom{+} \; \tfrac{1}{12} (a_{63} - a_{36} + a_{54} - a_{45}) e^1 +  \tfrac{1}{12} (b_{36} - b_{63} + b_{45} - b_{54}) e^7, \\
    \tau_2 =& \phantom{+} \; \, \tfrac{1}{3} (b_{45} - b_{54} + b_{36} - b_{63}) e^{12} +  \tfrac{1}{3} (a_{54} - a_{45} + a_{63} - a_{36}) e^{27}  \\
    &+  \tfrac{1}{3} (2 b_{36} - b_{63} + b_{54} - 2 b_{45}) e^{34} +  \tfrac{1}{3} (2 a_{36} + a_{63} - 2 a_{54} - a_{45}) e^{35}   \\
    &+  \tfrac{1}{3} (2 a_{63} + a_{36} - 2 a_{45} - a_{54}) e^{46} +  \tfrac{1}{3} (2 b_{63} + b_{36} + 2 b_{54} - b_{45}) e^{56},   \\
    \tau_3 =& \phantom{+} \; \, \tfrac{2}{7} \left(  c_{45} - c_{54} + c_{36} - c_{63} \right) e^{127} + \tfrac{1}{7} \left(5 c_{54} + 2 c_{45} - 5 c_{36} - 2 c_{63} \right) e^{135}    \\
        &+ \tfrac{1}{4} \left(b_{63} - b_{36} + b_{54} - b_{45} \right) e^{136} + \tfrac{1}{4} \left(b_{63} - b_{36} + b_{54} - b_{45} \right) e^{145}   \\
        &+ \tfrac{1}{7} \left(5 c_{45} + 2 c_{54} - 5 c_{63} - 2 c_{36} \right) e^{146} + \tfrac{1}{4} \left(3 a_{36} + 3 a_{45} + a_{63} + a_{54} \right) e^{234}   \\
        &+ \tfrac{1}{4} \left(3 b_{36} + b_{63} - 3 b_{54} - b_{45} \right) e^{235} + \tfrac{2}{7} \left(c_{54} - c_{45} + c_{63} - c_{36} \right) e^{236}   \\
        &+ \tfrac{2}{7} \left(c_{54} - c_{45} + c_{63} - c_{36} \right) e^{245} +  \tfrac{1}{4} \left(3 b_{63} + b_{36} - 3 b_{45} - b_{54} \right) e^{246}   \\
        &- \tfrac{1}{4} \left(3 a_{63} + a_{36} + 3 a_{54} + a_{45}  \right) e^{256} - \tfrac{1}{7} \left(5 c_{45} + 2 c_{54} + 5 c_{36} + 2 c_{63} \right) e^{347}   \\
        &+ \tfrac{1}{4} \left(a_{54} - a_{45} + a_{63} - a_{36} \right) e^{367} + \tfrac{1}{4} \left(a_{54} - a_{45} + a_{63} - a_{36} \right) e^{457}  \\ 
        &+ \tfrac{1}{7} \left(5 c_{54} + 2 c_{45} + 5 c_{63} + 2 c_{36} \right) e^{567},  
\end{align*}   
\end{proposition}

\indent Given that $\tau_3 \in \Span\{e^{127}, e^{135}, e^{136}, e^{145},  e^{146}, e^{234}, e^{235}, e^{236}, e^{245}, e^{246}, e^{256}, e^{347}, e^{367}, e^{457}, e^{567}\}$, equation \eqref{corolario cuentoso 3} from Corollary \ref{cor: cuentoso} entails the following result.

\begin{corollary} \label{cor: tau27 en el caso antidiagonal}
    The symmetric $2$-tensor $\tau_{27}$ of $(\mathfrak{g}_{A,B,C}, \varphi)$ in the antidiagonal case satisfies 
\begin{equation*}
    \tau_{27}(e_m, e_{9-m}) = 0 \quad \text{for all $3 \leq n \leq 6$}.
\end{equation*}    
\end{corollary}

\indent We apply Corollary \ref{cor: tau27 en el caso antidiagonal} to equation \eqref{eq: divergencia formula adiag} to establish the main result of this Section.

\begin{theorem} \label{thm: adiag_divergenceless}
    Let $A$, $B$, $C \in \mathfrak{sl}(4, \mathbb{R})$ be pairwise-commuting and antidiagonal. Then $(\mathfrak{g}_{A,B,C}, \varphi)$ is a $\mathrm{G}_2$-structure with divergence-free full torsion tensor.  
\end{theorem}  

\noindent {\it Acknowledgements.} I am very grateful with Dr. Adrián Andrada for his continued guidance during the preparation of this paper. I would also like to thank Dr. Jorge Lauret, for a suggestion of his originated this article. Special thanks go to Alejandro Tolcachier for thoughtful conversations and encouragement, and to Facundo Javier Gelatti for his patient and constant assistance in technical matters.

\appendix 

\section{Maple 2023 code} \label{section: maple 2023}

Formulas in Proposition \ref{prop: formulas N2}, \ref{prop: taus skew}, \ref{prop: taus diag}, and \ref{prop: taus adiag} where established (or verified) with help of Maple 2023 from Mapplesoft\textsuperscript{TM}, see https://www.maplesoft.com/products/Maple/. In this appendix we briefly describe the code used. \\

\indent We need the following packages:  
\begin{verbatim}
with(DifferentialGeometry);    with(LieAlgebras);   
with(liesymm);    with(Tensor);    with(LinearAlgebra); 
\end{verbatim}

\indent We begin by defining the following array:
\begin{verbatim}
C := Array(1 .. 7, 1 .. 7, 1 .. 7, 0);

C[3, 7, 3] := -a33; C[3, 7, 4] := -a43; C[3, 7, 5] := -a53; C[3, 7, 6] := -a63;
C[4, 7, 3] := -a34; C[4, 7, 4] := -a44; C[4, 7, 5] := -a54; C[4, 7, 6] := -a64;
C[5, 7, 3] := -a35; C[5, 7, 4] := -a45; C[5, 7, 5] := -a55; C[5, 7, 6] := -a65;
C[6, 7, 3] := -a36; C[6, 7, 4] := -a46; C[6, 7, 5] := -a56; C[6, 7, 6] := a33 + a44 + a55;

C[1, 3, 3] := b33; C[1, 3, 4] := b43; C[1, 3, 5] := b53; C[1, 3, 6] := b63;
C[1, 4, 3] := b34; C[1, 4, 4] := b44; C[1, 4, 5] := b54; C[1, 4, 6] := b64;
C[1, 5, 3] := b35; C[1, 5, 4] := b45; C[1, 5, 5] := b55; C[1, 5, 6] := b65;
C[1, 6, 3] := b36; C[1, 6, 4] := b46; C[1, 6, 5] := b56; C[1, 6, 6] := -b33 - b44 - b55;

C[2, 3, 3] := c33; C[2, 3, 4] := c43; C[2, 3, 5] := c53; C[2, 3, 6] := c63;
C[2, 4, 3] := c34; C[2, 4, 4] := c44; C[2, 4, 5] := c54; C[2, 4, 6] := c64;
C[2, 5, 3] := c35; C[2, 5, 4] := c45; C[2, 5, 5] := c55; C[2, 5, 6] := c65;
C[2, 6, 3] := c36; C[2, 6, 4] := c46; C[2, 6, 5] := c56; C[2, 6, 6] := -c33 - c44 - c55;
\end{verbatim}

\indent We can now use
\begin{verbatim}
    L1 := LieAlgebraData(C, Ex1);
\end{verbatim}
\noindent to define the Lie algebra $\mathfrak{g}_{A,B,C}$; it should be clear that the array C above describe the commutators from equation \eqref{eq: commutators g_A,B,C}. We now set
\begin{verbatim}
    DGsetup(L1);
\end{verbatim}

\indent We define the $\mathrm{G}_2$-structure $\varphi$ from equation \eqref{eq: G2-estructura} and the metric $g$ it induces explicitly: 
\begin{verbatim}
phi :=    ((theta1 &wedge theta2) &wedge theta7) + ((theta3 &wedge theta4) &wedge theta7) 
        + ((theta5 &wedge theta6) &wedge theta7) + ((theta1 &wedge theta3) &wedge theta5)
        - ((theta1 &wedge theta4) &wedge theta6) - ((theta2 &wedge theta3) &wedge theta6)
        - ((theta2 &wedge theta4) &wedge theta5);
g := CanonicalTensors("Metric", "bas", 7, 0);
\end{verbatim}

\indent The torsion forms $\tau_0$, $\tau_1$, $\tau_2$, $\tau_3$ can be computed from equations \eqref{eq: tau taus} as follows:
\begin{verbatim}
tau0 := 1/7*HodgeStar(g, ExteriorDerivative(phi) &wedge phi)
tau1 := -1/12*HodgeStar(g, HodgeStar(g, ExteriorDerivative(phi)) &wedge phi)
tau2 := 1/3*evalDG(3*(- HodgeStar(g, ExteriorDerivative(HodgeStar(g, phi))) 
                    + 4*HodgeStar(g, tau1 &wedge HodgeStar(g, phi))))
tau3 := evalDG(HodgeStar(g, ExteriorDerivative(phi)) - tau0*phi
            -3*HodgeStar(g, tau1 &wedge phi))
\end{verbatim}

\noindent These lead to equations \eqref{tau0 N2} to \eqref{tau3 N2}; as it was mentioned in the text, equations \eqref{tau0 N2}, \eqref{tau1 N2}, and \eqref{tau2 N2} are verified to coincide with \cite[Proposition 2.7]{N}. In order to obtain the results from Propositions \ref{prop: taus skew}, \ref{prop: taus diag}, or \ref{prop: taus adiag} we need to modify the array C, which can only be done manually, and re-run the previous commands. \\

\indent It is straightforward to verify Lemma \ref{lemma: cuentoso} and Corollary \ref{cor: cuentoso}, as they follow from just computations with wedge products. As for Corollary \ref{cor: tau27 en el caso general}, we set
\begin{verbatim}
    tau27 := (x, y) -> HodgeStar(g, (Hook(x, phi) &wedge Hook(y, phi)) &wedge tau3)
\end{verbatim}
\noindent and evaluate $\tau_{27}(x,y)$ with $x = e_k$, $y = e_i$ for all $k = 1, 2, 7$ and $3 \leq i \leq 6$ one by one. We can calculate the elements in equation \eqref{eq: los tau27 clave} just as easily and, with that, find an explicit formula for $\DV T$; however, as it is mentioned in the  text, it appears that little gain is made in doing so, not to mention that it is perhaps not pleasing to depend much on software.  

\end{document}